\documentclass[twoside,10pt,leqno]{amsart}
\everymath={\displaystyle}
\usepackage{amsfonts}
\usepackage{amsmath}
\usepackage{amscd}
\usepackage{amssymb}
\usepackage{amsthm}
\usepackage{amsrefs}
\begin{document}
\newtheorem{theorem}{Theorem}
\newtheorem{lemma}{Lemma}
\newtheorem{corollary}{Corollary}
\newtheorem{conjecture}{Conjecture}
\newtheorem{prop}{Proposition}
\numberwithin{equation}{section}
\newcommand{\dif}{\mathrm{d}}
\newcommand{\ord}{\mathrm{ord}}
\newcommand{\lcm}{\mathrm{lcm}}
\newcommand{\intz}{\mathbb{Z}}
\newcommand{\ratq}{\mathbb{Q}}
\newcommand{\natn}{\mathbb{N}}
\newcommand{\comc}{\mathbb{C}}
\newcommand{\rear}{\mathbb{R}}
\newcommand{\prip}{\mathbb{P}}
\newcommand{\uph}{\mathbb{H}}

\title[Large Sieve Inequality in High Dimensions]{\bf An Improvement of a Large Sieve Inequality in High Dimensions}
\author{Liangyi Zhao}
\date{\today}
\maketitle

\begin{abstract}
In this paper, we present an improvement of a large sieve type inequality in high dimensions and discuss its implications on a related problem.
\end{abstract}

\section{Introduction}

It was in 1941 that J. V. Linnik \cite{JVL1} first introduced the idea of large sieve in the investigation of the distribution of quadratic non-residues.  Applications of the idea abound. \newline

The large sieve inequality, the present form of which was first introduced by H. Davenport and H. Halberstam \cite{DH1}, is stated as follows.  There are many references on the subject.  See, for example, \cites{EB2, MBB1, HM2, HM, PXG}.  We shall henceforth refer to it as the classical large sieve inequality.  For notational convenience, a set of real numbers $\{ x_k \}$ is said to be $\delta$-spaced modulo 1 if $\| x_j-x_k \| > \delta$, for all $j \neq k$, where henceforth if $x = ( x_1 , \cdots , x_n ) \in \rear^n$ $\|x\|$ denotes $\max _i \min_{k\in \intz} |x_i-k|$. \newline

\begin{theorem}[Classical Large Sieve Ineqaulity] \label{classicls}
Let $\{ a_n \}$ be an arbitrary set of complex numbers, $\{ x_k \}$ be a set of real numbers that is $\delta$-spaced modulo 1, and $M \in \intz$, $N \in \natn$.  Then
\begin{equation} \label{classiclseq}
\sum_k \left| \sum_{n=M+1}^{M+N} a_n e(x_kn) \right|^2 \ll (\delta^{-1}+N) \sum_{n=M+1}^{M+N} |a_n|^2,
\end{equation}
where the implied constant is absolute.
\end{theorem}
Save for the more precise implied constant, the above inequality is the best possible.  Montgomery and Vaughan \cite{HLMRCV} showed that
\begin{equation} \label{montvaughan}
\sum_k \left| \sum_{n=M+1}^{M+N} a_n e \left( x_k n \right) \right|^2 \leq \left( \delta^{-1}+N \right) \sum_{n=M+1}^{M+N} |a_n|^2,
\end{equation}
while Paul Cohen and Selberg have shown independently that $\delta^{-1}+N$ can be replaced by $\delta^{-1}+N-1$ which is absolutely the best possible, since Bombieri and Davenport \cite{BD1} gave examples of $\{x_k\}$ and $a_n$, with $\delta \to 0$, $N \to \infty$ and $N \delta \to \infty$ such that equality holds in \eqref{montvaughan} with $\delta^{-1}+N-1$.  However, in our paper, we shall not be concerned with the implied constants. \newline

As corollaries to Theorem~\ref{classicls}, we have the following inequality for additive characters.
\begin{equation} \label{classiclseqadd}
\sum_{q=1}^Q \sum_{\substack{a \; \bmod{ \; q} \\ \gcd(a,q)=1}} \left| \sum_{n=M+1}^{M+N} a_n e \left( \frac{a}{q} n \right) \right|^2 \ll (Q^2+N) \sum_{n=M+1}^{M+N} |a_n|^2.
\end{equation}
Various extensions of these classical results restricted to various kinds of special characters are also known \cites{LZ1, LZ4}.  Results similar to \eqref{classiclseq} are also known in higher dimension and the proof is also similar.  The following is quoted from \cite{MNH2}.
\begin{theorem} \label{lshighdim}
Let
\[ S(x_1, \cdots, x_k) = \sum_{n_1, \cdots, n_k} c(n_1, \cdots, n_k) e(n_1x_1 + \cdots +n_kx_k), \]
where the summation is over integer points in $k$ dimensional rectangle.  $M_j < n_j \leq M_j+N_j$ for $j=1, \cdots, k$.  Let $x^{(1)}, \cdots, x^{(R)}$ be real $k$ dimensional vectors, say $x^{(r)}=(x_1^{(r)},\cdots,x_k^{(r)})$, which satisfy
\[ \max_j \delta^{-1}_j \| x_j^{(r)}-x_j^{(s)} \| >1, \]
for all $r \neq s$, where $\delta_j$'s are positive numbers not exceeding $\frac{1}{2}$.  Then we have
\begin{equation}\label{lshighdimeq}
\sum_{r=1}^R \left| S(x^{(r)}) \right|^2 \leq \prod_{j=1}^k \left( \sqrt{N_j} + \sqrt{\delta_j^{-1}} \right)^2 \sum_{n_1, \cdots, n_k} |c(n_1, \cdots, n_k)|^2.
\end{equation}
\end{theorem}
Note that the length of the outer summation on the left-hand side of \eqref{classiclseq} does not exceed $\delta^{-1}$ and analogous statements can be made about that of \eqref{lshighdimeq}, and that the right-hand sides of \eqref{classiclseq}, \eqref{classiclseqadd} and \eqref{lshighdimeq} are essentially the sum of the lengths of the summations on the left-hand sides times the square of the $l_2$-norm of the sequence $\{a_n\}$.  It is this feature of the classic theorem that motivated our desire for improving the following, which is quoted from P. X. Gallagher \cite{PXG3}.
\begin{theorem} \label{lshighdim1}
Let $\cdot$ denote the usual dot product in $\rear^n$ and $c(a)$ be a complex-valued function on $\intz^n$.  Then
\begin{equation} \label{lshighdimeq1}
\sum_{\substack{\beta \in \rear^n/\intz^n \\ \ord(\beta) \leq X}} \left| \sum_{\substack{\alpha=(\alpha_1, \cdots, \alpha_n)\in \intz^n \\ \max_{1\leq i \leq n}|\alpha_i|\leq N}} c(\alpha) e ( \alpha \cdot \beta) \right|^2 \ll \left( N^n + X^{2n} \right) \sum_{\substack{\alpha=(\alpha_1, \cdots, \alpha_n)\in \intz^n \\ \max_{1\leq i \leq n}|\alpha_i|\leq N}} |c(\alpha)|^2,
\end{equation}
where the implied constant depends on $n$.
\end{theorem}
Here and after, $\ord(\beta)$ denotes the additive order of $\beta$ in $\rear^n/\intz^n$.  Hence if $\beta = \left( \frac{a_1}{q_1}, \cdots, \frac{a_n}{q_n} \right)$, then $\ord (\beta) = \lcm (q_1, \cdots, q_n)$.  \eqref{lshighdimeq1} follows easily from \eqref{lshighdimeq}.  But the result of Theorem~\ref{lshighdim} is more general than that of Theorem~\ref{lshighdim1}, as the outer summation of \eqref{lshighdimeq} can be considerably longer than that of \eqref{lshighdimeq1}.  Hence it is believed that $X^{2n}$ on the right-hand side of \eqref{lshighdimeq1} can replaced, as noted before, by $X^{n+1}$, the length of the outer summation(see Lemma~\ref{spacing}), in spirit analogous to that of the classical large sieve inequality.  In other words, the majorant in \eqref{lshighdimeq1} might be replaced by
\begin{equation} \label{goal1}
\left( N^n + X^{n+1} \right) \sum_{\substack{\alpha=(\alpha_1, \cdots, \alpha_n)\in \intz^n \\ \max_{1\leq i \leq n}|\alpha_i|\leq N}} |c(\alpha)|^2.
\end{equation}
It is clear that both terms above are necessary.  Set $c(\alpha)=1$ for all $\alpha$ and $N=1$, we see that $X^{n+1}$ is needed.  Taking $X=1$ and $c(\alpha)=e(-\alpha \cdot \beta)$ gives the conclusion that $N^n$ is necessary.  However, \eqref{goal1} is not enough.  The following is a counter example.  Let
\begin{equation*}
T = \{ \beta \in \rear^n/\intz^n \; : \; \ord(\beta) \leq X, \beta = ( \beta_1, \cdots , \beta_n), \; \beta_i=0 \; \mbox{for} \; 2 \leq i \leq n, \; \beta=\frac{a_1}{q_1}, \; q_1 \in \prip \},
\end{equation*}
where $\prip$ denotes the set of prime numbers.  It is clear that $T$ is of size $\gg X^{2-\epsilon}$, as it can be identified with the Farey fractions of level $X$ with prime denominators.  Let $c(\alpha)=1$ for all $\alpha$.  We have
\begin{eqnarray}
\nonumber & \sum_{\substack{\beta \in \rear^n/\intz^n \\ \ord(\beta) \leq X}} \left| \sum_{\substack{\alpha=(\alpha_1, \cdots, \alpha_n)\in \intz^n \\ \max_{1\leq i \leq n}|\alpha_i|\leq N}} c(\alpha) e ( \alpha \cdot \beta) \right|^2 \\
\nonumber \geq & \sum_{\beta \in T} \left| \sum_{\substack{\alpha=(\alpha_1, \cdots, \alpha_n)\in \intz^n \\ \max_{1\leq i \leq n}|\alpha_i|\leq N}} c(\alpha) e ( \alpha \cdot \beta) \right|^2 \\
\nonumber = & N^{2n-2} \sum_{\substack{1\leq p \leq X \\ p \in \prip}} \sum_{\substack{a \mod p \\ \gcd(a,p)=1}} \left| \sum_{|m|\leq N} e \left( \frac{a}{q} m \right) \right|^2 \\
\nonumber = & N^{2n-2} \sum_m \sum_{m'} \sum_{\substack{1\leq p \leq X \\ p \in \prip}} \sum_{\substack{a \mod p \\ \gcd(a,p)=1}} e \left( \frac{a}{p} (m-m') \right) \\
\nonumber = & N^{2n-2} \left[ \sum_m \sum_{m'} \sum_{\substack{1\leq p \leq X \\ p \in \prip \\ p|(m-m')}} (p-1) -  \sum_m \sum_{m'} \sum_{\substack{1\leq p \leq X \\ p \in \prip \\ p\nmid (m-m')}} 1 \right] \\
\nonumber = & N^{2n-2} \left[ \sum_{\substack{1\leq p \leq X \\ p \in \prip}} \mathop{\sum_m \sum_{m'}}_{p|(m-m')} p - \sum_{\substack{1\leq p \leq X \\ p \in \prip}} \sum_m \sum_{m'} 1 \right] \\
\label{counterexample2} \geq & c(\epsilon) N^{2n-1} X^{2-\epsilon} - \pi(X)N^{2n},
\end{eqnarray}
for some $c(\epsilon)>0$ that depends only on $\epsilon$ and as usual $\pi(x)$ denotes the number of primes not exceeding $x$.  But \eqref{goal1} gives the majorant of
\[ N^{2n}+N^{n}X^{n+1}. \]
Taking $N=X^{1+\theta}$ for any $0 < \theta <1$, which ensures the dominance of the positive term in \eqref{counterexample2}, we see that the majorant of \eqref{goal1} is not enough. \newline

The following notations and conventions are used throughout paper. \newline

\noindent $e(z) = \exp (2 \pi i z) = e^{2 \pi i z}$. \newline
$f = O(g)$ means $|f| \leq cg$ for some unspecified postive constant $c$. \newline
$f \ll g$ means $f=O(g)$. \newline
$f \asymp g$ means $f \ll g$ and $g \ll f$.  Unless otherwise stated, all implied constants in $\ll$, $O$ and $\asymp$ are absolute. \newline
$\qed$ denotes the end of a proof or the proof is easy and standard.

\subsection*{Acknowledgment}The author wishes to thank Professors P. X. Gallagher and J. B. Friedlander, the former for suggesting the problem and both for the helpful discussions.  The author was supported by a grant from the Faculty Development and Research Fund at the United States Military Academy and a post-doctoral fellowship at the University of Toronto during this work.

\section{Preliminary Lemmas}

In this section, we quote the lemmas needed for the results of this paper.  As in the best-known proof of the classical large sieve inequality, we need the duality principle.

\begin{lemma}[Duality Principle]\label{dual}
Let $T=[t_{mn}]$ be a square matrix with entries from the complex numbers.  The following two statements are equivalent:
\begin{enumerate}
\item For any absolutely square summable sequence of complex numbers $\{ a_n
\}$, we have
\begin{equation} \label{dual1}
\sum_m \left| \sum_n a_n t_{mn} \right|^2 \leq D \sum_n |a_n|^2.
\end{equation}
\item For any absolutely square summable sequence of complex numbers
$\{ b_n \}$, we have
\begin{equation} \label{dual2}
\sum_n \left| \sum_m b_m t_{mn} \right|^2 \leq D \sum_m |b_m|^2.
\end{equation}
\end{enumerate}
\end{lemma}
\begin{proof}  This is a standard result.  See Theorem 288 in \cite{GHHJELGP}. \end{proof}

We shall also need the following lemma regarding the spacing of certain $n$ dimensional vectors.  Here and after, we set
\[ S= \left\{ \beta \in \rear^n/\intz^n : \ord ( \beta ) \leq X, \; \beta = \left(\frac{a_1}{q_1}, \; \cdots \; \frac{a_n}{q_n}\right), \; \frac{X}{2} \leq q_1 \leq X \right\}. \]

\begin{lemma} \label{spacing}
Let $\epsilon >0$ be given and $Y>0$
\[ M(X,Y) = \max_{\beta \in S} \# \left\{ \beta' \in S \; : \; \|\beta-\beta'\|<Y \right\}. \]
Then we have
\begin{equation} \label{spacingeq}
M(X,Y) \ll X^{\epsilon}\left(X^{n+1}Y^n + X^2Y + 1\right),
\end{equation}
where the implied constant depends on $n$ and $\epsilon$.
\end{lemma}
\begin{proof} We estimate the size of the set of our interest in the following way.  Fix $\beta = \left( \frac{a_1}{q_1}, \cdots, \frac{a_n}{q_n} \right) \in S$.  The number of $\frac{a_1'}{q_1'}$'s with $1 \leq a_1' < q_1'$, $X/2 < q_1' \leq X$ and $\gcd(a_1',q_1')=1$ such that $\left\| \frac{a_1}{q_1}-\frac{a_1'}{q_1'} \right\| <Y$ does not exceed $X^2Y+1$.  For each such $\frac{a_1'}{q_1'}$, we have the following number of choices for the other coordinates of $\beta'$.
\begin{equation*}
\sum_{i=1}^{[X/q_1']+1} \left( \sum_{k|iq_1'} (kY+1) \right)^{n-1} \ll \sum_{i=1}^{[X/q_1']+1} \left( Y^{n-1} (iq_1')^{n-1+\epsilon} + (iq_1')^{\epsilon} \right) \ll Y^{n-1} X^{n-1+\epsilon} + X^{\epsilon}.
\end{equation*}
Recall that $X/2 \leq q_1' \leq X$.  Hence in total, we have
\[ M(X,Y) \ll Y^nX^{n+1+\epsilon}+Y^{n-1}X^{n+\epsilon}+YX^{2+\epsilon}+X^{\epsilon}. \]
The term $Y^{n-1} X^{n+\epsilon}$ is not necessary, for $Y^{n-1} X^{n+\epsilon} \geq YX^{2+\epsilon}$ implies $XY\geq 1$ and hence $Y^{n-1} X^{n+\epsilon} \leq Y^nX^{n+1+\epsilon}$.  Hence the result follows.
\end{proof}
Note that upon taking $Y=1$, we get the the size of the set $S$ is $O_{\epsilon}(X^{n+1+\epsilon})$.  Therefore, in the light of Lemma~\ref{spacing}, so long as $Y$ is not so small that no regularity of distribution of elements of $S$ can be expected, the spacing property of $S$ is essentially as expected, as given in the first term of \eqref{spacingeq}. \newline

It is somewhat a melancholy admission, as will be noted in Section 4, that the term $X^{2+\epsilon}Y$ is necessary in Lemma~\ref{spacing}.  The following is an example to that effect.  Let 
\begin{equation*}
T'= \{ \beta \in \rear^n/\intz^n \; : \; \ord(\beta) \leq X, \beta = ( \beta_1, \cdots , \beta_n), \; \beta_i=0 \; \mbox{for} \; 2 \leq i \leq n \},
\end{equation*}
The spacing properties of elements in $T'$ are the same as those of the Farey fractions of level $X$.  Hence
\[ \max_{\beta \in T'} \# \left\{ \beta' \in T' \; : \; \| \beta - \beta' \| < Y \right\} \asymp X^2Y+1. \]
Therefore, we have $M(X,Y) \gg X^2Y$.  Taking $Y=X^{-1-\theta}$ for any $0 < \theta <1$, we see that the term $X^2Y$ is needed in \eqref{spacingeq}.

\section{Main Contention}

The objective is to have an upper bound for the following sum.
\[ \sum_{\substack{\beta \in \rear^n/\intz^n \\ \ord(\beta) \leq X}} \left| \sum_{\substack{\alpha=(\alpha_1, \cdots, \alpha_n)\in \intz^n \\ \max_{1\leq i \leq n}|\alpha_i|\leq N}} c(\alpha) e ( \alpha \cdot \beta) \right|^2. \]
Without trying too hard and in the light of Lemma~\ref{spacing}, simply applying Cauchy's inequality would give us the majorant of
\[ N^nX^{n+1+\epsilon} \sum_{\alpha} |c(\alpha)|^2, \]
which is already better than \eqref{lshighdimeq1} when $N^n \ll X^{n-1-\epsilon}$.  But certainly we hope to do better than {\it this}.  Furthermore, some applications require that the size of $X$ is well controlled.  To that end, we have the following.

\begin{theorem} \label{improve}
Under the notations that have been in use thus far, we have
\begin{equation} \label{improveeq}
\sum_{\substack{\beta \in \rear^n/\intz^n \\ \ord(\beta) \leq X}} \left| \sum_{\substack{\alpha=(\alpha_1, \cdots, \alpha_n)\in \intz^n \\ \max_{1\leq i \leq n}|\alpha_i|\leq N}} c(\alpha) e ( \alpha \cdot \beta) \right|^2 \ll X^{\epsilon} \left( X^{n+1}+N^{n-1}X^2 + N^n \right) \sum_{\alpha} |c(\alpha)|^2,
\end{equation}
where the implied constant depends on $n$ and $\epsilon$.
\end{theorem}
\begin{proof}
It will suffice to break up the outer sums into dyadic intervals.  Together with the application of the duality principle, Lemma~\ref{dual}, it suffices to show that
\begin{equation} \label{dualapply}
\sum_{\alpha} \left| \sum_{\beta \in S} b(\beta) e( \alpha \cdot \beta ) \right|^2 \ll X^{\epsilon} \left( X^{n+1}+N^{n-1}X^2 + N^n \right)  \sum_{\beta} |b(\beta)|^2,
\end{equation}
for any sequence of complex numbers $\{ b(\beta) \}$ and where $S$ and the summations over $\alpha$ and $\beta$ are as before.  \newline

Set $\phi(x) = \left( \frac{\sin \pi x}{2x} \right)^2$.  By positivity, the left-hand side of \eqref{dualapply} is bounded above by
\[ \sum_{\alpha \in \intz^n} \prod_{i=1}^n \phi \left( \frac{\alpha_i}{2N} \right) \left| \sum_{\beta\in S} b(\beta) e( \alpha \cdot \beta ) \right|^2, \]
where the sum over $\alpha$ is now extended over all elements of $\intz^n$.  Expanding the modulus square in the above and factoring, it becomes
\begin{eqnarray}
\nonumber & \sum_{\beta \in S} \sum_{\beta' \in S} b(\beta) \bar{b}(\beta') \sum_{\alpha \in \intz^n} \prod_{i=1}^n \phi \left( \frac{\alpha_i}{2N} \right) e( \alpha_i (\beta_i-\beta_i') ) \\
\label{est1} = & \sum_{\beta \in S} \sum_{\beta' \in S} b(\beta) \bar{b}(\beta') \prod_{i=1}^n \sum_{\alpha_i=-\infty}^{\infty} \phi \left( \frac{\alpha_i}{2N} \right) e( \alpha_i (\beta_i-\beta_i') )
\end{eqnarray}
Set $V(y) = \sum_{n=-\infty}^{\infty} \phi \left( \frac{n}{2N} \right) e(ny)$.  Recall that the Fourier transform of $\phi(x)$ is precisely $\Lambda(s)=\max(1-|s|,0)$.  Hence, we apply the Poisson summation formula and a change of variables to obtain
\begin{eqnarray*}
V(y) & = & 2N \sum_{m=-\infty}^{\infty} \Lambda(2N(m+y)) \\
 & = & \frac{\pi^2 N}{2} \sum_{|m+y|<(2N)^{-1}} (1-2N|m+y|) \\
 & = & \frac{\pi^2 N}{2} (1-2N\|y\|),
\end{eqnarray*}
if $\|y\|<(2N)^{-1}$ and $V(y)=0$ otherwise.  Therefore, \eqref{est1} is 
\begin{eqnarray*}
= & \left( \frac{\pi^2 N}{2}\right)^n \mathop{\sum_{\beta} \sum_{\beta'}}_{\|\beta-\beta'\|<(2N)^{-1}} b(\beta) \bar{b}(\beta') \prod_{i=1}^n (1-2N\| \beta_i - \beta_i'\|) \\
\leq & \left( \frac{\pi^2 N}{2}\right)^n \mathop{\sum_{\beta} \sum_{\beta'}}_{\|\beta-\beta'\|<(2N)^{-1}} | b(\beta) \bar{b}(\beta') | \\
\leq & \left( \frac{\pi^2 N}{2}\right)^n \sum_{\beta} |b(\beta)|^2 M(X,(2N)^{-1}),
\end{eqnarray*}
with $M(X,Y)$ defined as in Lemma~\ref{spacing}.  Upon inserting the result of Lemma~\ref{spacing} with $Y=(2N)^{-1}$ and summing up all the dyadic intervals for $X$, our contention follows.
\end{proof}
From the discussion and the examples given in section 1, we can infer that the inequality in \eqref{improveeq} is essentially the best possible.

\section{Notes}

It was the inequality \eqref{lshighdimeq1} that was the starting point for P. X. Gallagher \cite{PXG3} in improving an estimate on the number $E_n(N)$ of monic polynomials 
\[ F(x) = X^n + a_1X^{n-1}+ \cdots + a_n \]
with integer coefficients and of height, $H(F)= \max(|a_1|, \cdots, |a_n|)$ not exceeding $N$ for which the Galois group is a proper subgroup of the symmetric group.  The problem was first studied by van der Waerden \cite{BLVDW} and improvements were later made by Knobloch \cites{HWK1, HWK2}.  Gallagher's improvement gives the bound
\[ E_n(N) \ll N^{n-\frac{1}{2}} \log N, \]
with the implied constant depending on $n$.  The size of $X$ required to ensure the dominance of $N^n$ in \eqref{lshighdimeq1} is the key factor for determining the negative part of the exponent of $N$ above.  Unfortunately, our result \eqref{improveeq} requires the exact same size for $X$ to ensure the dominance of $N^n$ and hence it leads to essentially the same bounds for $E_n(N)$ as above.

\bibliography{biblio}

\begin{bibdiv}
\begin{biblist}

\bib{MBB1}{article}{
    author={Barban, M.~B.},
     title={The "large sieve" method and its applications in the theory of
  numbers},
      date={1966},
   journal={Uspehi Matematiks Nauk},
    volume={21},
     pages={51\ndash 102},
}

\bib{EB2}{article}{
    author={Bombieri, E.},
     title={Le grand crible dans la th\'eorie analytique des nombres},
      date={1974},
   journal={Soci\'et\'e Mathematics France},
    volume={18},
}

\bib{BD1}{article}{
    author={Bombieri, E.},
    author={Davenport, H.},
     title={Some inequalities involving trigonometrical polynomials},
      date={1969},
   journal={Annali Scuola Normale Superiore - Pisa},
    volume={23},
     pages={223\ndash 241},
}

\bib{DH1}{article}{
    author={Davenport, H.},
    author={Halberstam, H.},
     title={The values of a trigonometric polynomial at well spaced points},
      date={1966},
   journal={Mathematika},
    volume={13},
     pages={91\ndash 96},
      note={{\it Corrigendum and Addendum}, Mathematika {\bf 14} (1967),
  232-299},
}

\bib{PXG}{article}{
    author={Gallagher, P.~X.},
     title={The large sieve},
      date={1967},
   journal={Mathematika},
    volume={14},
     pages={14\ndash 20},
}

\bib{PXG3}{inproceedings}{
    author={Gallagher, P.~X.},
     title={The large sieve and probabilistic Galois theory},
organization={American Mathematical Society},
      date={1973},
 booktitle={Proceedings of symposium on pure mathematics},
    volume={XXIV},
     pages={91\ndash 101},
}

\bib{GHHJELGP}{book}{
    author={Hardy, G.~H.},
    author={Littlewood, J.~E.},
    author={P\'olya, G.},
     title={Inequalities},
 publisher={Cambridge University Press},
      date={1964},
}

\bib{MNH2}{article}{
    author={Huxley, M.~N.},
     title={The large sieve inequality for algebraic number fields},
      date={1968},
   journal={Mathematika},
    volume={15},
     pages={178\ndash 187},
}

\bib{HWK1}{article}{
    author={Knobloch, H.-W.},
     title={Zum hilbertschen irreduzibilit\"{a}tssatz},
      date={1955},
   journal={Abh. Math. Sem. Univ. Hamburg.},
    volume={19},
     pages={176\ndash 190},
}

\bib{HWK2}{article}{
    author={Knobloch, H.-W.},
     title={Die seltenheit der reduziblen polynome},
      date={1956},
   journal={Jber. Deutsch. Math. Verein},
    volume={59},
     pages={12\ndash 19},
}

\bib{JVL1}{article}{
    author={Linnik, J.~V.},
     title={The large sieve},
      date={1941},
   journal={Doklady Akademii Nauk Soiuza Sovetskikh Sotsialisticheskikh
  Respublik},
    volume={36},
     pages={119\ndash 120},
      note={(Russian)},
}

\bib{HM}{book}{
    author={Montgomery, H.~L.},
     title={Topics in Multiplicative Number Theory},
    series={Lecture Notes in Mathematics},
 publisher={Spring-Verlag},
   address={Barcelona, etc.},
      date={1971},
    volume={227},
}

\bib{HM2}{article}{
    author={Montgomery, H.~L.},
     title={The Analytic Priciples of Large Sieve},
      date={1978July},
   journal={Bulletin of the American Mathematical Society},
    volume={84},
    number={4},
     pages={547\ndash 567},
}

\bib{HLMRCV}{article}{
    author={Montgomery, H.~L.},
    author={Vaughan, R.~C.},
     title={The large sieve},
      date={1973},
   journal={Mathematika},
    volume={20},
     pages={119\ndash 134},
}

\bib{BLVDW}{article}{
    author={van~der Waerden, B.~L.},
     title={Die seltenheit der gleichungen mit affekt},
      date={1934},
   journal={Math. Ann.},
    volume={109},
     pages={13\ndash 16},
}

\bib{LZ1}{article}{
    author={Zhao, L.},
     title={Large sieve inequality for characters to square moduli},
      date={2004},
   journal={Acta Arithmetica},
    volume={112},
    number={3},
     pages={297\ndash 308},
}

\bib{LZ4}{article}{
    author={Zhao, L.},
     title={Large sieve inequality for special characters to prime square
  moduli},
      date={2004},
   journal={Functiones et Approximatio Commentarii Mathematici},
    volume={XXXII},
     pages={1\ndash 8},
}

\end{biblist}
\end{bibdiv}
\bibliographystyle{amsxport}

\noindent{Department of Mathematics} \newline
{University of Toronto} \newline
{100 Saint George Street} \newline
{Toronto, ON M5S 3G3 Canada}\newline
{\sc Email Address:} {\tt lzhao@math.toronto.edu} \newline

\end{document}